\documentclass[11pt]{article}

\usepackage{amsmath,amssymb,enumerate,pb-diagram,dsfont}
\usepackage[alphabetic]{amsrefs}
\def\today{\number\day .\number\month .\number\year}
\usepackage[applemac]{inputenc}

\def \1{\mathds{1}}
\def \al{\alpha}

\def \C{{\mathbb C}}

\def \df{\ \stackrel{\mbox{\rm\tiny def}}{=}\ }
\def \e{\emph}

\def \ga{\gamma}
\def \Ga{\Gamma}
\def \hom{{\rm hom}}
\def \H{\operatorname{H}}
\def \Hom{\operatorname{Hom}}

\def \om{\omega}

\def \ph{\varphi}
\def \prf{{\bf Proof: }}

\def \qed{\hfill $\square$}
\def \R{{\mathbb R}}

\def \res{\operatorname{res}}
\def \Sh{\operatorname{Sh}}

\def \Z{{\mathbb Z}}
\def \({\left(}
\def \){\right)}
\def \={{\ =\ }}

\newcommand{\tto}[1]{\stackrel{#1}{\longrightarrow}}

\renewcommand{\sp}[1]{\left\langle #1\right\rangle}
\newcommand{\ol}[1]{\overline{#1}}

\newtheorem{lemma}{Lemma}[section]

\newtheorem{exmples}[lemma]{Examples}

\newtheorem{exmple}[lemma]{Example}

\newtheorem{defi}[lemma]{Definition}

\newtheorem{theorem}[lemma]{Theorem}

\begin{document}

\pagestyle{myheadings} \markright{ITERATED INTEGRALS}

\title{Iterated Integrals and higher order invariants}
\author{Anton Deitmar \& Ivan Horozov\\ \ \\
Can. J. Math. 65, 544-552 (2013)}
\date{}
\maketitle

{\bf Abstract:}
We show that higher order invariants of smooth functions can be written as linear combinations of full invariants times iterated integrals.
The non-uniqueness of such a presentation is captured in the kernel of the ensuing map from the tensor product. This kernel is computed explicitly.
As a consequence, it turns out that higher order invariants are a free module of the algebra of full invariants.


$$ $$

\tableofcontents

\newpage
\section*{Introduction}

Modular forms of higher order have been studied extensively in recent years \cites{ChintaDiam,ES,HighOrd,DDDir,DKMO,DO,DS,DiaKleb,DiaSim,Sree}.
To construct them, one often uses iterated integrals. Form dimension formulae \cite{DiaSim} it is clear that in the case of holomorphic forms, iterated integrals do not give all higher order forms.
But there are strong indicators \cites{DKMO,DO,DS,Sree} that in the case of smooth functions, all higher order forms are indeed obtainable via iterated integrals.
This has been an implicit open question for a while which is answered affirmatively in the present paper.
It is shown that on any smooth manifold, the smooth module of higher order invariants is generated by the space of homotopy invariant iterated integrals, thus allowing to deduce structure assertions on higher order invariants from iterated integrals.

In Section 1 we recall the ingredients of the theory of iterated integrals needed in the sequel.
In Section 2 we show that the restriction map from free homotopy invariant iterated integrals to loops is a surjective map.
This assertion is the key to the further sections.
In Section 3 we first show that homotopy invariant iterated integrals always give higher order invariants.
For the case of surfaces, this was proven by Sreekantan in \cite{Sree}.
We then use the result from Section 3 to finally deduce the main result, saying that higher order invariants can alway be expressed by iterated integrals.
In the final section we consider restricted forms of higher order, which come with triviality assumptions along boundary components. 
This is typical for automorphic forms, where these restrictions refer to cusp forms or cuspidal cohomology.
We formulate the corresponding theorem asserting that the main result is stable under boundary restrictions.

\section{Generalities on iterated integrals}
In this section we fix notations.
Let $X$ be a smooth connected manifold and $x_0,x\in X$ points.
We write $PX$ for the \e{path space}, i.e., the set of all smooth maps $p:[0,1]\to X$.
We also write $PX_{x_0}$ for the subset of all paths that start at $x_0$ and $PX_{x_0,x}$ for the subset of all smooth paths from $x_0$ to $x$.
The space $LX_{x_0}=PX_{x_0,x_0}$ is also called the \e{loop space} at $x_0$.

For a path $p$ and $1$-forms $\om_1,\dots,\om_r$ we define the \e{iterated integral:}
$$
\int_p\om_1\cdots\om_r=\int_0^1\int_0^{t_r}\dots\int_0^{t_2}p^*\om_1(t_1)\, p^*\om_2(t_2)\dots p^*\om_r(t_r).
$$
For an integer $s$, let $B_s(X)$ denote the space of all maps $\om: PX\to\C$ which are linear combinations of iterated integrals of length $\le s$.
Here we include constants as they may be considered as iterated integrals of length zero.
We also write $B(X)$ for the union of all $B_s(X)$ as $s$ varies.
Let
$$
T(\Omega^1(X))=\C\oplus \Omega^1(X)\oplus\left[\Omega^1(X)\otimes \Omega^1(X)\right]\oplus\dots
$$
be the tensorial algebra over the space $\Omega^1(X)$ of smooth $1$-forms. 
The map assigning $\om_1\otimes\dots\otimes\om_r$ to the map $p\mapsto\int_p\om_1\cdots\om_r$ is a linear map from $T(\Omega^1(X))$ to $B(X)$.
This map has a non-trivial kernel which has been determined by Chen in \cite{Alg}.

We denote by $B_s(X)_{x_0}$ the set of restrictions of elements of $B_s(X)$ to $PX_{x_0}$ and the space $B_s(X)_{x_0,x}$ is defined analogously.

\begin{lemma}\label{lem1.1}
\begin{enumerate}[\rm (a)]
\item If $\ph$ is an orientation preserving diffeomorphism on $[0,1]$, then $\int_p\om_1\cdots\om_r=\int_{p\circ\ph}\om_1\cdots\om_r$.
\item If $F$ is  a diffeomorphism on $X$, then $\int_{F\circ p}\om_1\dots\om_r=\int_p(F^*\om_1)\dots (F^*\om_r)$.
\item If $p$ and $q$ are composable paths, then 
$$
\int_{pq}\om_1\cdots\om_r=\sum_{j=0}^r\int_p \om_1\cdots\om_j\int_q\om_{j+1}\cdots\om_r.
$$
\item
One has
$$
\(\int_p\om_1\cdots\om_r\)\(\int_p\om_{r+1}\cdots\om_{r+s}\)
=\sum_\sigma\int_p\om_{\sigma(1)}\cdots\om_{\sigma(r+s)},
$$
where the sum runs over all $(r,s)$-shuffles, i.e., permutations $\sigma$ on $r+s$ letters with
$\sigma^{-1}(1)<\dots<\sigma^{-1}(r)$ and $\sigma^{-1}(r+1)<\dots<\sigma^{-1}(r+s)$.
\item $\int_{p^{-1}}\om_1\cdots\om_r=(-1)^r\int_p\om_r\cdots\om_1$, where $p^{-1}(t)=p(1-t)$.
\item For given $\om\in B(X)$, we extend the map $p\mapsto\int_p\om$ to the free abelian group $\Z[PX]$ generated by $PX$.
For given $\al_1,\dots,\al_s\in PX_{x_0,x_0}$ let $\eta=(\al_1-1)(\al_2-1)\dots (\al_s-1)\in\Z[PX]$.
For 1-forms $\om_1,\dots,\om_r$, $r\le s$ we have
$$
\int_\eta\om_1\cdots\om_r=\begin{cases}
\prod_{i=1}^s\int_{\al_i}\om_i& r=s,\\
0&r<s.
\end{cases}
$$
\end{enumerate}
\end{lemma}

Note that (d) implies that
$$
B(X)_{x_0}\=\bigcup_{s=0}^\infty B_s(X)_{x_0}
$$
is a filtered algebra and likewise for $B(X)_{x_0,x}$.

\prf (a)-(e) are easy exercises.
(f) is a result of Hain, see \cite{Hain}, Proposition 2.13.
\qed

If we replace the tensor product on $T(X)$ by the \e{shuffle product} $*$ given by
$$
\om_1\cdots\om_r*\om_{r+1}\cdots\om_{r+s}
=\sum_\sigma\om_{\sigma(1)}\cdots\om_{\sigma(r+s)},
$$
where the sum runs over all $(r,s)$-shuffles, we obtain the \e{shuffle algebra} $\Sh(X)$.
We have shown that the iterated integrals form an algebra homomorphism
$$
\Sh(X)\to B(X),
$$
where the latter is an algebra under pointwise multiplication.

Let $B_s(X)^\hom$ denote the space of all elements of $B_s(X)$ which are invariant under homotopies with fixed end-points.
Similarly define $B_{s}(X)^\hom_{x_0}$ and $B_{s}(X)^\hom_{x_0,x}$.

\section{The restriction map}

\begin{theorem}
The restriction map $B_s(X)_{x_0}^\hom\to B_s(X)_{x_0,x_0}^\hom$ is surjective.
\end{theorem}

\prf
For a real vector space $V$, let $T(V)=\R\oplus V\oplus V^{\otimes 2}\oplus\dots$ be the tensorial algebra.
Consider the map
$D:T(\Omega^1(X))\to T(\Omega(X))$ given by
\begin{eqnarray*}
D(\om_1\otimes\dots\otimes\om_n)
&=&
\sum_{j=1}^n\om_1\otimes\dots\otimes (d\om_j)\otimes\dots\otimes\om_n\\
&&
+\sum_{j=1}^{n-1}\om_1\otimes\dots\otimes(\om_j\wedge\om_{j+1})\otimes\dots\om_n,
\end{eqnarray*}
and $D(c)=0$, for $c\in\R$.
By Proposition 1.5.2 of \cite{ChenIPI} the map $\int:T(\Omega^1X)\to \Omega^0(PX)$ sending $\om_1\otimes\dots\otimes\om_s$ to the map $p\mapsto\int_p\om_1\cdots\om_s$ extends to a map $T(\Omega X)\to \Omega(PX)$ such that for $\om\in T(\Omega^1(C))$ one has
$$
d\(\int\om\)=-\int D\om-p_0^*\om_1\int\om_1\cdots\om_s+(-1)^{s-1}p_1^*\om_s\int\om_1\cdots\om_{s-1},
$$
where $p_0,p_1:PX\to X$ are the evaluation maps that map a path to its start and end point, respectively.
Therefore, the kernel of $D$ is mapped to homotopy invariant iterated integrals on $PX_{x_0}$.
Let $M$ be the kernel of $D$.
Then $M$ has a natural filtration by degrees:
$$
\R=M_0\subset M_1\subset\dots
$$
The map $I_{x_0}:\om_1\otimes\dots\otimes\om_s\mapsto \int\om_1\dots\om_s$ maps $M_s$ to $B_s(X)_{x_0}^\hom$.
We get a commutative diagram:
$$
\begin{diagram}
\node{M_s}\arrow{e,t}{I_{x_0}}\arrow{se,b}{I_{x_0,x_0}}	
	\node{B_s(X)_{x_0}^\hom}\arrow{s,r}{\rm res}\\
\node[2]{B_s(X)_{x_0,x_0}^\hom.}
\end{diagram}
$$
The integral $\int\om_1\dots\om_s$ is a function on the path space $PX$.
Its restriction to the loop space $L_{x_0}X$ is homotopy invariant, if and only if it is locally constant, which is the case if and only if it is anihilated by the differential of the complex $\Lambda(PM)$ as in \cite{ChenIPI}.
Now the differential $D$ above also coincides with the differential of the bar construction on \cite{ChenIPI}, Section 4.1.
 Theorem 4.1.1 of \cite{ChenIPI} states that the iterated integral map is an isomorphism of graded differential algebras from that bar construction to the iterated integrals on the loop space.
Therefore the iterated integral map $I_{x_0,x_0}$ is surjective, hence the restriction map is surjective, too.
\qed

\section{The fundamental group}

For a group $\Ga$ we write its group ring as $A=\Z\Ga$. Let $J\subset\Z\Ga$ be the augmentation ideal, i.e., the span of all elements of the form $(\ga -1)$, where $\ga\in\Ga$.
For any $\Z\Ga$-module $V$ we write $\H_s^0(\Ga,V)$ for the $\Z$-module of all $v\in V$ with $J^sv=0$.
This space can be identified with $\Hom_{\Z\Ga}(\Z\Ga/J^s,V)$.
The elements of $\H_s^0(\Ga,V)$ for varying $s$ are called \e{higher order invariants}.
If $v$ is in $\H_s^0(\Ga,V)$, but not in $\H_{s-1}^0(\Ga,V)$, then $s$ is called the \e{order} of $v$.

Let $X$ be a connected smooth manifold, $x_0\in X$ a base-point, and $\Ga=\pi_1(X,x_0)$ the corresponding fundamental group.
We consider $\Ga$ as group of deck transformations on the universal covering $\tilde X$ of $X$.
We also fix a pre-image $x_0$ in $\tilde X$, which we will denote by the same symbol $x_0$ as no confusion can arise.

As $\tilde X$ is simply connected, the iterated integral $\int_p\om$ for $\om\in B_s(\tilde X)^\hom$ only depends on the endpoints $x,y$ of the path $p$.
We therefore write $\int_{x}^{y}\om=\int_p\om$.

Every $\ga\in\Ga$ can be viewed as a homotopy class of a loop based at $x_0\in X$.
In this way we get a map $B_s(X)_{x_0,x_0}^\hom\to {\rm Map}(\Ga,\C)$ that maps $\om\in B_s(X)_{x_0,x_0}^\hom$ to the map $\ga\mapsto\int_\ga\om$.
The latter map induces a $\Z$-linear map from the group ring $\Z\Ga$ to $\C$.
It is the content of Chen's de Rham Theorem for fundamental groups (see \cite{ChenIPI}, Corollar 1 to Theorem 2.6.1, see also \cite{PS}) that this map induces a bijection
$$
B_s(X)_{x_0,x_0}^\hom\tto\cong\Hom_\Z(\Z\Ga/J^{s+1},\C).
$$

Each $\om\in B_s(X)_{x_0}^\hom$ lifts to $\tilde X$ and gives an element of $B_s(\tilde X)_{x_0}^\hom$.
For $x\in\tilde X$ we write $\int_{x_0}^x\om$ for the iterated integral of this lift over any path joining $x_0$ and $x$, or, what amounts to the same, for the integral of $\om$ over the projection to $X$ of any such path.

\begin{theorem}\label{thm2.1}
If $\om\in B_s(X)_{x_0}^\hom$, then the function $\int_{x_0}^x\om$, $x\in\tilde X$, is an invariant of order at most $s+1$ in the $\Ga$-module $C^\infty(\tilde X)$.
This defines an injective linear map
$$
\Psi:B_s(X)_{x_0}^\hom\hookrightarrow \H_{s+1}^0(\Ga, C^\infty(\tilde X)).
$$
\end{theorem}

The case when $X$ is the hyperbolic plane is in the paper \cite{Sree}.

\prf
Let $\om\in B_s(X)^\hom$ and for $x\in\tilde X$ set $f_\om(x)=\int_{x_0}^x\om$.
We have to show
$$
\left[(\ga_1-1)\cdots(\ga_{s+1}-1)\right]^*f_\om =0
$$
for any $\ga_1,\dots,\ga_{s+1}\in\Ga$.
For given $x\in X$ and $\ga\in\Ga$, we choose a path 
$\ga_x$ from $x$ to $\ga x$.
The map $\ga\mapsto \ga_x$ is extended linearly to a map $\Z\Ga\to\Z[PX]$.
For every $x\in X$ we also fix a smooth path $p_x$ from $x_0$ to $x$.
Let $\om\in B_s(X)^\hom$ and let $\eta=\sum_\ga c_\ga\ga$ be an arbitrary element of the group ring $\Z\Ga$.
We have
\begin{eqnarray*}
\eta^*f_\om(x) &=& \sum_\ga c_\ga\ga^*f_\om(x)= \sum_\ga c_\ga\int_{x_0}^{\ga x}\om\\
&=& \sum_\ga c_\ga\int_{p_x\ga_x}\om= \int_{p_x\sum_\ga c_\ga\ga_x}\om=\int_{p_x\eta_x}\om.
\end{eqnarray*}
We apply this to the element $(\ga_1-1)\cdots (\ga_{s+1}-1)$ of the group ring and we look at any monomial $\om_1\cdots\om_r$ in $\om$, where the $\om_j$ are $1$-forms on $X$.
We then have
$$
\int_{p_x[(\ga_1-1)\dots(\ga_{s+1}-1)]_x}\om_1\dots\om_r=\sum_{k=0}^r\int_{p_x}\om_1\dots\om_k\int_{[(\ga_{1}-1)\dots(\ga_{s+1}-1)]_x}\om_{k+1}\dots\om_r.
$$
Let $\bar x\in X$ be the image of $x\in\tilde X$ and let $\ga_{\bar x}$ be the image of $\ga_x$ in $X$.
Then $\ga_{\bar x}$ is a loop based at $\bar x$.
As the forms $\om_j$ are $\Ga$-invariant, we have
\begin{eqnarray*}
\int_{[(\ga_{1}-1)\dots(\ga_{s+1}-1)]_x}\om_{k+1}\dots\om_r
&=&
\int_{[(\ga_{1}-1)\dots(\ga_{s+1}-1)]_{\bar x}}\om_{k+1}\dots\om_r\\
&=&\int_{(\ga_{1,\bar x}-1)\dots(\ga_{s+1,\bar x}-1)}\om_{k+1}\dots\om_r\\
&=& 0
\end{eqnarray*}
by Lemma \ref{lem1.1} (e).
This proves the first claim.
For the injectivity of the induced map let $\om\in B_s(X)_{x_0}^\hom$ with $\int_{x_0}^x\om=0$.
This just means that $\om=0$ in $B_s(X)_{x_0}^\hom$.
\qed

$ $


We formally set $H_0^0=0$ and $B_{-1}=0$.

\begin{theorem}\label{thm3.2}
Let $K$ be the ideal in the algebra $B(X)_{x_0}$ generated by $d(C^\infty(X))$ and let $K_s^\hom=B_s(X)_{x_0}^\hom\cap K$.
Then  $K_s^\hom$ is the kernel of the restriction map
$$
B_s(X)_{x_0}^\hom\tto{\res}B_s(X)_{x_0,x_0}^\hom,
$$
and
$\Psi$ induces an isomorphism of $C^\infty(X)$-modules
$$
C^\infty(X)\otimes\(B_s(X)_{x_0}^\hom/K_s^\hom\)\ \tto\cong\ H_{s+1}^0(\Ga,C^\infty(\tilde X)).
$$
It follows that $H_{s+1}^0(\Ga,C^\infty(\tilde X))$ is a free module of the algebra $C^\infty(X)$ of smooth functions.
\end{theorem}

\prf
The fact that $K_s^\hom$ is the kernel of the restriction map is a consequence of Theorem 4.5 in \cite{Alg}.
Write $\bar B_s=B_s/B_{s-1}$ and let
$$
\bar K_s^\hom=\ker\left[B_s(X)_{x_0}^\hom\to\bar B_s(X)_{x_0,x_0}^\hom\right].
$$
Then $\bar K_s^\hom$ contains $B_{s-1}(X)_{x_0}^\hom$ and our assertion is equivalent to
$$
C^\infty(X)\otimes\(B_s(X)_{x_0}^\hom/\bar K_s^\hom\)\ \tto\cong\ \bar \H_{s+1}^0(\Ga,C^\infty(\tilde X)),
$$
which is what we prove.

By Chen's de Rham Theorem for fundamental groups (see \cite{ChenIPI}, Corollar 1 to Theorem 2.6.1, see also \cite{PS}) the evaluation of iterated integrals gives an isomorphism
$$
\bar B_s(X)_{x_0,x_0}^\hom\tto\cong \Hom_{\Z}(J^s/J^{s+1},\C).
$$
The right hand side can also be viewed as $\Hom_{\Z\Ga}(J^s/J^{s+1},\C)$ and as such be embedded into
$$
\Hom_{A}(J^s/J^{s+1},C^\infty(X))\cong\Hom_{A}(J^s/J^{s+1},C^\infty(\tilde X)),
$$
where we have written $A=\Z\Ga$.
More precisely, the image in 
$$
\Hom_{A}(J^s/J^{s+1},C^\infty(X))\cong C^\infty(X)\otimes \Hom_{A}(J^s/J^{s+1},\C)
$$
is a basis of this $C^\infty(X)$-module, which means that we have an isomorphism of $C^\infty(X)$-modules,
$$
C^\infty(X)\otimes \bar B_s(X)_{x_0,x_0}^\hom\tto\cong
\Hom_{A}(J^s/J^{s+1},C^\infty(\tilde X)).
$$

\begin{lemma}
The cohomology group $H^1(\Ga,C^\infty(\tilde X))$ is trivial.
\end{lemma}

\prf
A 1-cocycle is a map $\al:\Ga\to C^\infty(\tilde X)$ such that $\al(\ga\tau)=\ga\al(\tau)+\al(\ga)$ holds for all $\ga,\tau\in\Ga$.
We have to show that for any given such map $\al$ there exists $f\in C^\infty(\tilde X)$ such that $\al(\tau)=\tau f-f$.

Fix a smooth map $u:\tilde X\to [0,1]$ such that
$$
\sum_{\tau\in\Ga}u(\tau^{-1} x)\ \equiv\ 1,
$$
where we can assume that the sum is locally finite.
Set
$$
f(x)\= -\sum_{\tau\in\Ga}\al(\tau)(x)\,u(\tau^{-1} x).
$$
Then the function $f$ lies in the space $C^\infty(\tilde X)$.
We now compute for $\ga\in\Ga$,
\begin{eqnarray*}
\ga f(x)- f(x) &=& f(\ga^{-1} x)-f(x)\\
&=& \sum_{\tau\in\Ga}\al(\tau)(x) u(\tau^{-1} x) -
\al(\tau)(\ga^{-1} x)u(\tau^{-1}\ga^{-1} x)\\
&=& \sum_{\tau\in\Ga}\al(\tau) (x) u(\tau^{-1}x) 
+\al(\ga)(x)\sum_{\tau\in\Ga} u((\ga\tau)^{-1}x)\\
&& -\sum_{\tau\in\Ga}\al(\ga \tau)(x)u((\ga\tau)^{-1} x)
\end{eqnarray*}
The first and the last sum cancel and the middle sum is $\al(\ga)(x)$.
Therefore, the  lemma is proven.
\qed

Since $H^1(\Ga,C^\infty(\tilde X))=0$,  Lemma 2.1 in \cite{CompactQuot} implies that the exact sequence
$$
0\to J^s/J^{s+1}\to A/J^{s+1}\to A/J^s\to 0,
$$
induces an isomorphism
\begin{eqnarray*}
\Hom_{A}(J^s/J^{s+1},C^\infty(\tilde X))
&\cong& 
\underbrace{\Hom_{A}(A/J^{s+1},C^\infty(\tilde X))
/\Hom_{A}(A/J^{s},C^\infty(\tilde X))}_{\df\ol{\Hom}_{A}(A/J^{s+1},C^\infty(\tilde X))}\\
&\cong& \bar \H_{s+1}^0(\Ga,C^\infty(\tilde X)).
\end{eqnarray*}
We have to show that the ensuing diagram
$$
\begin{diagram}
\node{\bar B_s(X)_{x_0}^\hom}
		\arrow{e,t}{\res}
		\arrow{s}
	\node{\bar B_s(X)_{x_0,x_0}^\hom}
		\arrow{s,J}\\
\node{\ol{\Hom}_{A}(A/J^{s+1},C^\infty(\tilde X))}
		\arrow{e,t}{\cong}
	\node{\Hom_A(J^s/J^{s+1},C^\infty(\tilde X))}
\end{diagram}
$$
commutes.
This will give the claim, as we already have seen that the right vertical arrow becomes an isomorphism after tensoring with $C^\infty(X)$.
This commutativity is a direct consequence of formula (f) in Lemma \ref{lem1.1}.
\qed

\section{Parabolic restrictions}

To satisfy needs in the theory of automorphic forms, we also introduce further structure: Let $P\subset\Ga$ be a conjugation-invariant subset and let $\sp P$ be the normal subgroup generated by $P$.
Then $AJ_P$ is a two-sided ideal of $A$, where $J_P\subset \Z\sp P\subset A$ is the augmentation ideal of $\sp P$.
Let
$$
J_s\= J^s+AJ_P,
$$
and let
$$
\H_{P,s}^0(\Ga,V)\=\{ v\in V: J_s v=0\}.
$$
In the theory of automorphic forms, see \cite{ES, HighOrd}, $P$ will be the set of parabolic elements.

The $P$-restriction translates on the side of iterated integrals to the following.
Recall Chen's map
$$
B_s(X)_{x_0,x_0}^\hom\tto\cong\Hom_\Z(\Z\Ga/J^{s+1},\C).
$$
We define the space $B_{P,s}(X)_{x_0,x_0}^\hom$ to be the inverse image of 
$$
\Hom_\Z(\Z\Ga/J_{s+1},\C)\=\{\al\in\Hom_\Z(\Z\Ga/J^{s+1},\C): \al(p-1)=0\ \forall_{p\in P}\}
$$ 
under this map.
Finally, we set $B_{P,s}(X)_{x_0}^\hom$ to be the inverse image of $B_{P,s}(X)_{x_0,x_0}^\hom$ under the restriction map.
So then $B_{P,s}(X)_{x_0}^\hom$ is the set of all $\om\in B_{s}(X)_{x_0}^\hom$ with $\int_{p-1}\om=0$ for every $p\in P$.
Theorems \ref{thm2.1} and \ref{thm3.2} generalize to the following.

\begin{theorem}
If $\om\in B_{P,s}(X)_{x_0}^\hom$, then the function $\int_{x_0}^x\om$ is an invariant of order at most $s+1$ in the $\Ga$-module $C^\infty(\tilde X)$.
This defines an injective linear map
$$
\Psi:B_{P,s}(X)_{x_0}^\hom\hookrightarrow \H_{P,s+1}^0(\Ga, C^\infty(\tilde X)).
$$
\end{theorem}

\begin{theorem}
Let $K$ be the ideal in the algebra $B(X)_{x_0}$ generated by $d(C^\infty(X))$ and let $K_{P,s}^\hom=B_{P,s}(X)_{x_0}^\hom\cap K$.
Then  $K_{P,s}^\hom$ is the kernel of the restriction map
$$
B_{P,s}(X)_{x_0}^\hom\tto{\res}B_{P,s}(X)_{x_0,x_0}^\hom,
$$
and
$\Psi$ induces an isomorphism of $C^\infty(X)$-modules
$$
C^\infty(X)\otimes\(B_{P,s}(X)_{x_0}^\hom/K_{P,s}^\hom\)\ \tto\cong\ H_{P,s+1}^0(\Ga,C^\infty(\tilde X)).
$$
\end{theorem}
 \prf The proofs are the same as in the unrestricted case. One only has to check that the restriction conditions match, which is easy to see.
 \qed

{\small Mathematisches Institut,
Auf der Morgenstelle 10,
72076 T\"ubingen,
Germany,
\tt deitmar@uni-tuebingen.de\\ ivan.horozov@uni-tuebingen.de}

\today

\end{document}